\documentclass[12pt]{article} 
\usepackage{amstex, amsthm}

\newcommand{\doublespace}{
   \renewcommand{\baselinestretch}{1.2}
   \large\normalsize}

\renewcommand{\Bbb}{\mathbb}
\renewcommand{\frak}{\mathfrak}

\def \Z{\Bbb Z}
\def \C{\Bbb C}

\def \Q{\Bbb Q}

\def \wt{{\rm wt}}

\def \Res{{\rm Res}}

\def \End{{\rm End}}

\def \mod{{\rm mod}}
\def \<{\langle} 
\def \>{\rangle} 

\def \a{\alpha }

\def \l{\lambda }
\def \L{\Lambda }

\def \b{\beta }

\def \o{\omega}

\newcommand{\1}{\mathbf{1}}
\newcommand{\NO}{\,{\raise0.25em\hbox{$\mathop{\hphantom{\cdot}}%
\limits^{_{\circ}}_{^{\circ}}$}}\,}% The definition of the normal ordering.

\theoremstyle{plain}
 \newtheorem{theorem}{Theorem}[section]
	
	\newtheorem{lemma}[theorem]{Lemma}
	\newtheorem{proposition}[theorem]{Proposition}

 \newtheorem{remark}[theorem]{Remark}

\doublespace

\begin{document}

\newtheorem{thmm}{Theorem}
\newtheorem{thm}{Theorem}[section]
\newtheorem{prop}[thm]{Proposition}
\newtheorem{cor}[thm]{Corollary}
\newtheorem{lem}[thm]{Lemma}
\newtheorem{rem}[thm]{Remark}
\newtheorem{de}[thmm]{Definition}

\begin{center}
{\Large {\bf Representations of vertex operator algebra $V_L^+$ \\
for rank one lattice $L$}} \\
\vspace{0.5cm}
Chongying Dong\footnote{Supported by NSF grant 
DMS-9700923 and a research grant from the Committee on Research, UC Santa Cruz.}
and Kiyokazu Nagatomo\footnote{On leave of absence from
Department of Mathematics, Graduate School of Science,
Osaka University, Toyonaka, Osaka 560-0048, Japan.
This work was partly supported by Grant-in-Aid for Scientific Research,
the Ministry of Education, Science and Culture.}
\\
Department of Mathematics, University of California\\
 Santa Cruz, CA 95064\\
\end{center}

\hspace{1.5 cm}

\begin{abstract}
We classify the irreducible modules for 
the fixed point vertex operator subalgebra $V_L^+$ 
of the vertex operator algebra
$V_L$ associated to a positive definite even lattice 
of rank 1 under the automorphism lifted from the $-1$ isometry
of $L.$ 
\end{abstract}

\section{Introduction}\label{section:1}
\setcounter{equation}{0}

Vertex operator algebras $V_L$ associated to an arbitrary even positive
definite lattice $L$ have been studied extensively, and their representation
theory and fusion rules have been understood very well (see [B], [FLM2],
[D1], [DL1], [DLM1]). It is well known that 
the vertex operator algebra $V_L$ has an order $2$ automorphism $\theta$ 
which is deduced from the $-1$ isometry of the lattice [FLM2].
The $\theta$-invariants $V_L^+$ is a simple vertex operator subalgebra of
 $V_L.$ In this
paper we classify the irreducible modules for $V_L^+$ for all rank 1
lattice $L.$ The classification result says that any irreducible
module for $V_L^+$ is isomorphic to either a submodule of an irreducible
$V_L$-module or a submodule of an irreducible $\theta$-twisted
$V_L$-module. This confirms a conjecture in the orbifold conformal field 
theory [DVVV] in this special case.

The study of $V_L^+$ was initiated in [FLM1] during the course of
constructing the moonshine module although the notion of vertex operator
algebra was not available back then. It became clear later in [FLM2] (also
see [B]) that $V_{\L}^+$ (where $\L$ is the Leech lattice) 
is a vertex operator subalgebra of the
moonshine (module) vertex operator algebra $V^{\natural}$ which
is a direct sum of $V_{\L}^+$  and an irreducible module for $V_{\L}^+.$ 
In fact the moonshine module is the first example of so-called
``orbifold conformal field theory.'' The important role of 
the twisted modules in the orbifold theory was also noticed and
used in the construction of the moonshine module. 

There is a fundamental difference between vertex operator
algebras $V_L$ and $V_L^+.$ In order to see this we need to 
recall their construction. 
Set ${\frak h}=\C\otimes_{\Z}L$ and
$\hat{\frak h}={\frak h}\otimes \C[t,t^{-1}]\oplus \C c.$ Then the affine
Lie algebra $\hat{\frak h}$ has an automorphism $\theta$ of order 2 such that
$\theta (h\otimes t^n)=-h\otimes t^n$ and $\theta(c)=c.$ 
Then $V_L$ as a vector space is a tensor product of $M(1)$
and $\C[L]$ where $M(1)=S(
{\frak h}\otimes t^{-1}\C[t^{-1}])$ and $\C[L]$ is
the group algebra of $L.$ The map $\theta$ extends to an algebra automorphism
of  $S( {\frak h}\otimes t^{-1}\C[t^{-1}])$ and maps $e^{\alpha}$ to
$e^{-\alpha}$ for $\alpha\in L$ where we use $e^{\alpha}$ to denote
the corresponding element in $\C[L].$ It is clear that $V_L$ has the subspace
${\frak h}\otimes t^{-1}.$ So the affine
algebra $\hat{\frak h}$ with $c=1$ is a substructure of $V_L$ and
plays a very important role in the classification
of irreducible (twisted modules) for $V_L$ (see [D1] and [D2]). But
this affine algebra is not available to $V_L^+.$ This explains why
the representation theory for $V_L^+$ is more difficult.

Although $V_L^+$ provides a large class of concrete and important
examples of vertex operator algebras, the study of representation theory for 
an arbitrary $V_L^+$ is very limited
so far. It was proved in [DGH] that if $L$ contains a sublattice
of type $D_1^d$ then $V_L^+$ is rational but the classification of
irreducible modules even in this case remains open except when 
$L$ is the Leech lattice (see [D3]).  

The classification of irreducible modules for $V_L^+$ is well motivated
by the problem of classification of rational vertex operator 
algebras. The classification of rational vertex operator algebras 
is definitely one of the most important problems in the theory
of vertex operator algebra and it has immediate applications to
the classification of rational conformal field theory. 
If the central charge is less than 1, the classification
problem is not difficult as the vertex operator subalgebra generated
by the Virasoro element has only finitely many irreducible
modules. So the first nontrivial case is the central charge 
equal to 1. It is believed that
all rational vertex operator algebras of central charge 1 are
$V_L, V_L^+$ and $V_{L_2}^G$ where the rank of $L$ is 1, 
$L_2$ is the root lattice of type $A_1,$
$G$ is a finite subgroup of $SO(3)$ of type $E$ and $V_{L_2}^G$ 
is the corresponding invariants.  As we mentioned already that
the representation theory of $V_L$ including the fusion rules
is clear, one has to understand $V_L^+$ and $V_{L_2}^G$ better
and eventually characterize them. 

Another importance of studying the vertex operator algebra $V_L^+$ 
lies in the connection of $V_L^+$ with the $W$-algebra $W(2,4,k)$
(cf. [BFKNRV]) where the even positive integer $k$ is the half of
the square length of
generators of the lattice $L.$ It was pointed out in [DG] that 
$V_L^+$ is generated by the Virasoro vector and two additional
highest weight vectors for the Virasoro algebra of weights
$4$ and $k.$ So $V_L^+$ is the vertex operator algebra associated to
the $W$-algebra $W(2,4,k)$ of central charge $1.$ We may expect some
application of the results in this paper to the study of the $W$-algebra
$W(2,4,k).$  

The present paper is a continuation of [DN] in which we determine
the Zhu's algebra $A(M(1)^+)$ in the case $d=1$ and
classified the irreducible modules for $M(1)^+$ where
$M(1)^+$ is the $\theta$-invariants of $M(1).$
It is well known that $M(1)^+$ is a vertex operator subalgebra of
$V_L^+.$ The ideals, techniques and results in [DN] have been
extensively used and significantly extended in the present paper.
As in [DN] our main strategy is to determine the Zhu's
algebra $A(V_L^+)$ whose inequivalent irreducible modules
have a one to one correspondence with the inequivalent irreducible
(admissible) modules for $V_L^+.$ It turns out that $A(V_L^+)$ is
a semisimple commutative associative 
algebra of dimension $k+7$ generated by the image
of the three generators of $V_L^+$ in $A(V_L^+).$ 

The organization of this paper is as follows: In Section 2 we review 
the Zhu's algebra and related results. We also 
briefly review the vertex operator algebras $V_L,$ $V_L^+$ and
their (twisted) modules. In Section 3 we introduce
the three generators of $V_L^+$ following [DG] and 
give the commutator relations of the component operators of these
generators. We then show that how to obtain a ``small'' spanning 
set of $A(V_L^+).$ In section 4 we use a PBW type generating 
result to give an even ``smaller'' spanning set of $A(V_L^+)$ in terms
of the images of the three generators of $V_L^+.$ Section 5
is the core of this paper. In this section we first 
find the four relations among the three generators of $A(V_L^+)$ two
of which were in [DN] already. These relations
are good enough for us to determine a basis of $A(V_L^+)$
in the three separate cases: $k$ is not a perfect square;
$k$ is an even perfect square and $k$ is an odd perfect square different 
from 1. The case $k=1$ needs a special treatment although it is easy:
in this case $V_L^+$ is isomorphic to another lattice vertex
operator algebra $V_{L'}$ corresponding $k=4.$

\section{Preliminaries}\label{section:2}
\setcounter{equation}{0}

In this section after recalling
 the definitions of admissible modules for a vertex
operator algebra and a rational vertex operator algebra from [DLM2] and
[Z] we review the definition of Zhu's algebra and related results. 
We then review the construction of vertex operator algebra
$V_L$ associated to an even positive definite lattice of rank 1 and
its representations (see [B], [FLM2], [D1]). We also define the
automorphism $\theta$ of $V_L$ which is a lifting from
the $-1$ isometry of $L$ and the fixed point vertex operator subalgebra
$V_L^+.$  The $\theta$-twisted modules for $V_L$ are also discussed.

\subsection{Modules and Zhu's algebras}

We begin with a vertex operator algebra $V$ (cf. [B], [FLM2]) and
an automorphism of $g$ of $V$ of finite order
$T.$  Then $V$ decomposes into eigenspaces with respect to
the action of $g$ as $V=\bigoplus_{r\in \Z/T\Z}V^r$
where $V^r=\{v\in V|gv=e^{-2\pi ir/T}v\}$.

An admissible $g$-twisted $V$-module (cf. [DLM2], [Z]) is a $\frac{1}{T}\Z$-graded vector space
$$M=\sum_{n=0}^{\infty}M(\frac{n}{T}),$$
with top level 
$M(0)\ne 0,$ equipped 
with a linear map
$$\begin{array}{l}
V\longrightarrow (\End\,M)\{z\}\\
v\longmapsto\displaystyle{ Y_M(v,z)=\sum_{n\in\Q}v_nz^{-n-1}\ \ \ (v_n\in
\End\,M)}
\end{array}$$
satisfying the following conditions;
for all $0\leq r\leq T-1,$ $u\in V^r$, $v\in V,$ 
$w\in M$,
\begin{eqnarray*}
& &Y_M(u,z)=\sum_{n\in \frac{r}{T}+\Z}u_nz^{-n-1}, \label{1/2}\\ 
& &u_nw=0\ \ \                                  
\mbox{for}\ \ \ n\gg 0,\label{vlw0}\\
& &Y_M({\bf 1},z)=1,\label{vacuum}
\end{eqnarray*}
\[
\begin{array}{c}
\displaystyle{z^{-1}_0\delta\left(\frac{z_1-z_2}{z_0}\right)
Y_M(u,z_1)Y_M(v,z_2)-z^{-1}_0\delta\left(\frac{z_2-z_1}{-z_0}\right)
Y_M(v,z_2)Y_M(u,z_1)}\\
\displaystyle{=z_2^{-1}\left(\frac{z_1-z_0}{z_2}\right)^{-r/T}
\delta\left(\frac{z_1-z_0}{z_2}\right)
Y_M(Y(u,z_0)v,z_2)},
\end{array}
\]
where $\delta(z)=\sum_{n\in\Z}z^n$ and
all binomial expressions are to be expanded in nonnegative
integral powers of the second variable;
$$u_mM(n)\subset M(\wt(u)-m-1+n)$$
if $u$ is homogeneous. 
If $g=1$, this reduces to the definition of an admissible $V$-module.

A $g$-{\em twisted $V$-module} is
an admissible $g$-twisted $V$-module $M$ which carries a 
$\C$-grading  by weight. That is, we have
$$M=\coprod_{\lambda \in{\C}}M_{\lambda} $$
where $M_{\l}=\{w\in M|L(0)w=\l w\}.$ Moreover we require that 
$\dim M_{\l}$ is finite and for fixed $\l,$ $M_{{n\over T}+\l}=0$
for all small enough integers $n.$ Again if $g=1$  we get an ordinary
$V$-module.

A vertex operator algebra is called rational if any admissible module
is a direct sum of irreducible admissible modules. It was proved in
[Z] and [DLM2] that if $V$ is rational then $V$ has only finitely
many irreducible admissible modules and each irreducible
admissible module is an ordinary module. 

Zhu introduced an associative algebra $A(V)$ associated to a vertex
operator algebra $V$ 
which is extremely useful in study the representation theory of $V$ [Z]. 
In fact we will compute the Zhu's algebra $A(V_L^+)$ to determine
the irreducible modules for $V_L^+.$ 

Let $V$ be a vertex operator algebra. 
For homogeneous $u,v\in V,$ we define  products $u*v$ and $u\circ v$
 as follows:
\begin{eqnarray}\label{eqn:2.1}
& &u*v={\rm Res}_{z}\left(\frac{(1+z)^{{\rm wt}(u)}}{z}Y(u,z)v\right)
 =\sum_{i=0}^{\infty}{{\rm wt}(u)\choose i}u_{i-1}v\nonumber  \\
& &u\circ v={\rm Res}_{z}\left(\frac{(1+z)^{{\rm wt}(u)}}{z^2}Y(u,z)v
\right)
 =\sum_{i=0}^{\infty}{{\rm wt}(u)\choose i}u_{i-2}v. 
\end{eqnarray} 
Then extends (\ref{eqn:2.1}) to linear products on $V.$ Let $O(V)$ be
the linear span of $u\circ v$ for $u,v\in V.$ Set $A(V)=V/O(V).$

Let $M$ be an admissible module for $V.$ Following [DLM2] we define
the ``vacuum space'' 
$$\Omega(M)=\{w\in M|u_nw=0, u\in V, n\geq \wt(u)\}.$$
Then $\Omega(M)$ contains $M(0)$ and each $o(u)=u_{\wt(u)-1}$
for homogeneous $u\in V$ preserves $\Omega(M).$ One can extend
$o(u)$ to all $u\in V$ be the linearity.  
Then we have (see [Z], [DLM2])

\begin{theorem}\label{theorem:2.3}
{\rm (1)} $A(V)$ is an associative algebra under
multiplication $*$ and with identity $\1+O(V)$ and central element
$\omega+O(V).$ 

\noindent
{\rm (2)} The map $u\mapsto o(u)$ gives a representation of $A(V)$ on 
$\Omega(M)$
for any admissible $V$-module $M.$ 
Moreover, if $V$ is rational $A(V)$ is a finite dimensional semisimple algebra.

\noindent
{\rm (3)} The map  $M\to M(0)$ gives a bijection between the set of equivalence
classes of irreducible admissible $V$-modules and the  set of equivalence
classes of simple $A(V)$-modules.
\end{theorem}

Following [DN]  we write $[u]=u+O(V)\in A(V).$ 
We define $u\sim v$ for $u,v\in V$ if $[u]=[v].$ This induces a relation on 
$\End\, V$ such that for $f,g\in \End V$, $f\sim g$ if and only if
$fu\sim gu$ for all $u\in V$. 

We also need the following results from [Z].
\begin{proposition}\label{proposition:2.4}
{\rm (1)} Assume that $u\in V$ homogeneous, $v\in V$ and $ n \geq 0$.
Then 
\begin{equation}\label{eqn:2.2}
\Res_z\left(
\frac{(1+z)^{\wt(u)}}{z^{2+n}}Y(u,z)v
\right)
= \sum_{i=1}^\infty \binom{\wt(u)}{i}u_{i-n-2}v
\in O(V).
\end{equation}

\noindent
{\rm (2)} If $u$ and $v$ are homogeneous elements of $V$, then
\begin{equation}\label{eqn:2.3}
u*v\sim
\Res_z\left(
\frac{(1+z)^{\wt(v)-1}}{z}Y(v,z)u
\right).
\end{equation}

\noindent
{\rm (3)} For any $n\geq 1$, 
\begin{equation}\label{eqn:2.4}
L(-n)\sim (-1)^n
\left\{
(n-1)(L(-2)+L(-1))+L(0)
\right\}
\end{equation}
where $L(n)$ are the Virasoro operators given by 
$Y(\o,z)=\sum_{n\in\Z}L(n)z^{-n-2}.$

\noindent
{\rm (4)} For any $u \in V$, 
\begin{equation}\label{eqn:2.5}
[u]*[\omega] = [(L(-2)+L(-1))u].
\end{equation}
\end{proposition}

\subsection{Vertex operator algebras $V_L$ and $V_L^+$}
\label{subsection:2.1}

We work in the setting of [FLM2]. Let $L$ be an 
even lattice of rank one with nondegenerate symmetric ${\Bbb Z}$-bilinear 
form $\langle\cdot,\cdot\rangle$, ${\frak h}=L\otimes_{\Bbb Z}{\Bbb C}$
and $\hat{\frak h}_{\Bbb Z}$  the corresponding Heisenberg algebra.
Let $M(1)$ be the associated irreducible induced module for 
$\hat{\frak h}_{\Bbb Z}$ such that the canonical central element of 
$\hat{\frak h}_{\Bbb Z}$ acts as 1. Define $V_L=M(1)\otimes \C[L]$ where
$\C[L]$ is the group algebra of $L$ with a basis $\{e^{\a}|\a\in L\}.$
 Set $\1=1\otimes 1$ and $\omega=\frac{1}{2}\beta(-1)^2$ where $\beta\in\frak h$ such that 
$\<\b,\b\>=1.$  
It was proved in~[B] and [FLM2] that there is a linear map 
$$\begin{array}{lcr}
V_{L}&\longrightarrow& (\mbox{End}\,V_{L})[[z,z^{-1}]],\hspace*{3.6 cm} \\
v&\longmapsto& Y(v,z)=\displaystyle{\sum_{n\in\Z}v_nz^{-n-1}\ \ \ (v_n\in
\mbox{End}\,V_{L})}
\end{array}$$ such that
$V_{L}=(V_{L},Y,{\bf 1},\omega)$ is a simple vertex operator algebra.
Let $L^{\circ}=\{x\in\frak h\mid \<x,L\>\subset \Z\}$ be the dual lattice
of $L$. 
Then the irreducible modules for $V_L$ are $V_{L+\gamma}=M(1)\otimes
\C[L+\l]$ where $\l$ runs over the coset representatives
of $L$ in $L^\circ$ (see [D1]).  Moreover,
$V_L$ is a rational vertex operator algebra (see [DLM1]). To be more 
precise, let $L=\Z\a$ such that $\<\a,\a\>=2k.$ Then $L^\circ=\frac{1}{2k}L$
and the irreducible modules for $V_L$ are $V_{L+\frac{i}{2k}\alpha}$
for $i=0,...,2k-1.$ 

Let $\theta$ be the linear automorphism of $V_{L^\circ}$
such that $\theta (u\otimes e^{\gamma})
=\theta(u)\otimes e^{-\gamma}$ for $u\in M(1)$ and $\gamma\in L^\circ$.
Here the action of $\theta$ on $M(1)$ is given by $\theta(\a_1(n_1)\cdots
\a_k(n_k))=(-1)^k\a_1(n_1)\cdots\a_k(n_k)$. Then the restriction
of $\theta$ to $V_L$ is a VOA automorphism.  Let $M$ be an $\theta$-stable
subspace of $V_{L^\circ}.$ We denote the $\pm 1$ eigenspaces by $M^{\pm}$
respectively. Then $M(1)^+$ is a vertex operator subalgebra of $V_L^+.$

We have (see Theorems 4.4 and 6.1 of [DM])  
\begin{proposition}\label{proposition:2.1}
{\rm (1)} $V_L^+$ is a simple vertex operator algebra. 

\noindent
{\rm (2)} $V_{L}^{\pm}$ and $V_{L+\frac{1}{2}\a}^{\pm}$ are irreducible 
$V^+_L$-modules.

\noindent
{\rm (3)} $V_{L+\frac{i}{2k}}$ and $V_{L+\frac{2k-i}{2k}\alpha}$ are 
isomorphic and irreducible $V_L^+$-modules for $i=1,\dots,k-1.$ 
\end{proposition}

Next we discuss the $\theta$-twisted modules of $V_L$ following Chapter
9 of [FLM2]. Then $L/2L$ is an abelian group isomorphic to $\Z_2$ 
has two irreducible modules $T_1,T_2$ such that $\a+2L$ acts as scalars
$1$ and $-1$ respectively. Let $\hat {\frak h}[-1]$ be the twisted Heisenberg
algebra. As in Section 1.7 of [FLM2] we also denote by $M(1)$
the unique irreducible $\hat{\frak h}[-1]$-module with
the canonical central element acting by $1.$
Define the twisted space $V_{L}^{T_i}=M(1)\otimes T_i$.  It was
shown in [FLM2] and [DL2] that there is a linear map
\begin{eqnarray*}
V_{L}&\to& (\mbox{End}\,V_{L}^{T_i})[[z^{1/2},z^{-1/2}]],\hspace*{3.6 cm} \\
v&\mapsto& Y(v,z)=\displaystyle{\sum_{n\in\frac{1}{2}\Z}v_nz^{-n-1}}\ \
\end{eqnarray*}
such that $V_{L}^{T_i}$ is an irreducible
$\theta$-twisted module for $V_L$. Moreover, $V_{L}^{T_i}$
for $i=1,2$ give all irreducible $\theta$-twisted $V_L$-module
(see [D2]).

We also define a linear operator $\theta$ on $V_L^{T_i}$ such that
$$
\theta(\alpha_1(-n_1)\cdots\a_s(-n_s)\otimes t)=(-1)^{s}
\alpha_1(-n_1)\cdots\a_s(-n_s)\otimes t
$$
for $\a_i\in \frak h$, $n_i\in\frac{1}{2}+\Z$ and $t\in T_i$.
Then $\hat\theta_d Y(u,z)(\hat \theta_d)^{-1}=Y(\theta u,z)$ for $u\in V_L$
(cf.~[FLM2]).  
We have the decomposition $V_L^{T_i}=
(V_L^{T_i})^+\oplus(V_L^{T_i})^-$ where $(V_L^{T_i})^{\pm}$ 
are the $\pm 1$ eigenspaces 

Then we have (see [FLM2] and  Theorem 5.5 of [DLi])
\begin{proposition}\label{proposition:2.2} 
$(V_L^{T_i})^{\pm }$  are irreducible 
$V_L^+$-modules for $i=1,2.$
\end{proposition}

Our main result in this paper is that Propositions \ref{proposition:2.1} and 
\ref{proposition:2.2} give a complete list of irreducible modules for $V_L^+.$

\section{A spanning set of $A(V_L^+)$}
\label{section:3}
\setcounter{equation}{0}

In this section we use the ideas and techniques 
developed in [DN] to reduce the spanning set of
$A(V_L^+)$  to the images of $M(1)^+$ and $V_L^+(1)$ in $A(V_L^+).$
We shall use the vertex operators $Y(u,z)$ for $u\in V_L$ freely
and we refer the reader to [FLM2] for the definition of these
operators.  In Subsection 3.1 we review the bracket relations
for the component operators of the generators of $V_L^+.$ Subsection
3.2 gives several lemmas which are used in the later subsections.
In Subsections 3.3 and 3.4 
 we prove that the subspace $V_L^+(m)+O(V_L^+)$ of $A(V_L^+)$
can be generated by the subspace $M(1)^++V_L^+(1)+O(V_L^+)$ for 
all $m\geq 1.$ 

\subsection{The generators of $V_L^+$}
\label{subsetion:3.1new}

Recall from [DG] that the vertex operator algebra $M(1)^+$
is generated by $\o$ and
\begin{equation}
J = \b(-1)^4\1 -2\b(-3)\b(-1)\1 + \frac{3}{2}\b(-2)^2\1\label{eqn:3.3new}
\end{equation}
which is a singular vector of weight $4$ for the Virasoro algebra.
Also recall that $Y(\o,z)=\sum_{n\in\Z}L(n)z^{-n-2}$ and $J(z)=
\sum_{n\in\Z}J_nz^{-n-1}.$  The following lemma can be found in [DN].
\begin{lemma}\label{lemma:3.1new}
(1) For any $m,n\in\Z$,  
$$[L(m),J_n] = (3(m+1)-n)J_{n+m}.$$

\noindent(2) The commutators $[J_m,J_n]$ are expressed 
as linear combinations of
\[
L(p_1)\cdots L(p_s),\quad
L(q_1)\cdots L(q_t)J_r
\]
where $p_1,\dots,p_s,q_1,\dots,q_t,r\in\Z$ and $s,t\leq 3$.
\end{lemma}

For convenience we set 
\[
V_L^+(m)= M(1)^+\otimes (e^{m\alpha}+e^{-m\alpha})+
M(1)^-\otimes (e^{m\alpha}-e^{-m\alpha})
\]
for $m\geq 0.$ Then $ V_L^+(m) $ is an irreducible $M(1)^+$-module
and is also completely reducible module for the Virasoro algebra
\begin{equation}\label{eqn:3.2new}
V_L^+(m) =
\begin{cases}
\bigoplus_{p\geq 0}L(1,(m\sqrt{k}+p)^2)&\text{if $k\ne 0$ is a perfet square},\\
\bigoplus_{p\geq 0}L(1,(4p)^2)&\text{if $k=0$},\\
L(1,km^2) &\text{otherwise}
\end{cases}
\end{equation}
(cf. [DG]) where $L(1,h)$ is the highest weight module for the Virasoro
algebra with central charge 1 and highest weight $h.$ 

Set $E=e^{\a}+e^{-\a}.$ Then $V_L^+$ is generated by $\o,J$ and $E$  
(see Theorem 2.9 of [DG]). Since $E$ is a highest weight vector for
the Virasoro algebra of weight $k,$ we immediately have
$$[L(m),E_n] = ((k-1)(m+1)-n)E_{n+m}.$$ 
The commutator $[J_m,E_n]$ could be computed if one knows $J_sE$ for
$s\geq 0.$ Next we make a rough estimation of $J_sE.$ 
Since $\wt(J_sE) = k+3-s\leq k+3$ if $s\geq 0$, then $J_sE\in L(1,k)$ by 
(\ref{eqn:3.2new}) for $k>1.$ The following lemma now is obvious.
\begin{lemma}\label{lemma:3.2new}
Assume that $k>1.$ For nonnegative integer $n$, $J_nE$ is expressed as a linear
combinations of the set
\[
\left\{
L(-m_1)\cdots L(-m_s)E\,|\, m_1\geq m_2\geq \dots\geq m_1\geq 1,\,s\leq 3
\right\}.
\]
\end{lemma}

\subsection{Several lemmas}
\label{subsection:3.1}	

Recall from [FLM2] that $\alpha(z)=\sum_{m\in\Z}\alpha(m)z^{-m-1}$
for $\alpha\in\frak h$ and $Y(\alpha(-n-1)\1,z) =
\partial^{(n)}\alpha(z)$ for $n\in\Z_{\geq 0}$ where 
$\partial^{(n)}=\frac{1}{n!}\frac{d}{dz}.$ Then 
\begin{equation}\label{eqn:3.1}
\partial^{(n)}\alpha(z)
= \sum_{j\geq 0}\binom{-j-1}{n}\alpha(j)z^{-j-n-1}
+\sum_{j\leq -n-1}\binom{-j-1}{n}\alpha(j)z^{-j-n-1}.
\end{equation}
Set 
\[
E^m = e^{m\alpha}+e^{-m\alpha},\quad F^m =  e^{m\alpha}-e^{-m\alpha}
\]
for any integer $m$. Then $E^1$ is the $E$ defined in Subsection 3.1.
Also set $F=F^1.$ Notice that
\[
\alpha(0)E^m =2kmF^m,\quad \alpha(0)F^m = 2kmE^m.
\]

\begin{lemma}\label{lemma:3.1}
For $m,n\geq 1$,
\begin{align*}
(\alpha(-n)\alpha&(-1)\1)*E^m\\
&= \alpha(-n)\alpha(-1)E^m+2mk(n+(-1)^{n+1})\alpha(-n-1)F^m
+\sum_{i=0}^n c_i\alpha(-i)F^m
\end{align*}
for some $c_i\in \C$.
\end{lemma}
\begin{proof}
Let $v = \alpha(-n)\alpha(-1)\1.$ By (\ref{eqn:3.1}) we have
\begin{equation}\label{eqn:3.2}
Y(v,z) = \sum_{
\begin{Sb}
j_1,\,j_2\in\Z\\
\text{$j_1\geq 0$ or $j_1\leq -n$}
\end{Sb}}
\binom{-j_1-1}{n-1}
\NO \alpha(j_1)\alpha(j_2)\NO
z^{-j_1-j_2-n-1}.
\end{equation}
Recall that $Y(v,z)=\sum_{j\in\Z}v_jz^{-j-1}.$ Then
\[
v_{-1} = \sum_{
\begin{Sb}
j_1+j_2 = -n-1\\
\text{$j_1\geq 0$ or $j_1\leq -n$}
\end{Sb}}
\binom{-j_1-1}{n-1}
\NO\alpha(j_1)\alpha(j_2)\NO.
\]
Hence
\begin{equation*}
v_{-1}E^m=\left(
(-1)^{n-1}\alpha(-n-1)\alpha(0)+n\alpha(-n-1)\alpha(0)
+\alpha(-n)\alpha(-1)
\right)
E^m.
\end{equation*}
This proves that
\[
v_{-1}E^m
= \alpha(-n)\alpha(-1)E^m
+2mk(n+(-1)^{n+1})\alpha(-n-1)F^m.
\]

Next we consider $v_iE^m$ for $0\leq i\leq n$. By (\ref{eqn:3.1}), we see
\[
v_i =\sum_{
\begin{Sb}
j_1+j_2 =i-n\\
\text{$j_1\geq 0$ or $j_1\leq -n$}
\end{Sb}}
\binom{-j_1-1}{n-1} 
\NO \alpha(j_1)\alpha(j_2)\NO.
\]
So
\begin{align*}
v_iE^m &= \left(
(-1)^{n+1}\alpha(i-n)\alpha(0)
+\delta_{i0}\alpha(-n)\alpha(0)
\right)E^m\\
&=2km\left(
(-1)^{n+1}\alpha(i-n)+\delta_{i0}\alpha(-n)
\right)
F^m.
\end{align*}
Since $\wt(v) = n+1$ we see from (\ref{eqn:2.1}) that 
$$v*E^m
= \sum_{i=0}^{n+1}\binom{n+1}{i}v_{i-1}E^m.$$
Substitute the explicit expressions of $v_{i-1}$ into the equation
above to get the desired result.
\end{proof}

We say an element 
$u = \alpha(-n_1)\cdots\alpha(-n_r)v$  ($n_i>0,$ $v=E^m$ or $F^m$)
has the length $r$ with respect to $\alpha$ and we write
$\ell_\alpha(u) = r$. In general if $u$ is a linear combination
of such vectors $u^i$'s we define
the length of $u$ to be the maximal length among $\ell_{\a}(u^i).$

\begin{lemma}\label{lemma:3.2}
Let $n_1,n_2,\dots,n_r\in \Z_{>0}$ with $r$ even. Then
\[
(\alpha(-n_1)\cdots\alpha(-n_r)\1)*E^m
 =\alpha(-n_1)\cdots\alpha(-n_r)E^m + u
\]
where $u\in V_L^+(m)$ and $\ell_\alpha(u)<r$.
\end{lemma}
\begin{proof}
Let $v = \alpha(-n_1)\cdots\alpha(-n_r)\1$. 
By the definition of a vertex operator and (\ref{eqn:3.1}), we have
\[
Y(v,z) = \sum_{m_i\in\Z}c_{m_1n_1}\cdots c_{m_rn_r}
\NO \alpha(m_1)\cdots \alpha(m_r)\NO z^{-m-n}
\]
where $m= m_1+\cdots+m_r, n=n_1+\cdots+n_r$ and $c_{mn} = \binom{-m-1}{n-1}$,
$m_i\geq 0$ or $m_i\leq -n_i$. 
Therefore for $j\geq 0$ 
\[
v_{j-1}E^m
= \sum_{
\begin{Sb}
m= j-n\\
\text{$m_i\geq 0$ or $m_i\leq -n_i$}
\end{Sb}
}c_{m_1n_1}\cdots c_{m_rn_r}
\NO \alpha(m_1)\cdots \alpha(m_r)\NO E^m.
\]
If $j=0$ then either $m_i = -n_i$ for all $i$ or there exits $i$ such that
$m_i\geq 0$ . So in this case 
\[
v_{-1}E^m = \alpha(-n_1)\cdots \alpha(-n_r)E^m
+ u
\]
where $\ell_\alpha(u)<r$. If $j>0$ then there exists $i$ such that
 $m_i\geq 0$. 
This implies $\ell_\alpha(v_{j-1}E^m)<r.$
The lemma follows from the definition of $*$ product. 
\end{proof}

A similar argument gives:

\begin{lemma}\label{lemma:3.3} Let $n_1,...,n_r\in \Z_{>0}$ with
$r$ odd. Then 
\[
(\alpha(-n_1)\cdots\alpha(-n_{r-1})\1)*(\alpha(-n_r)F^m)
=\alpha(-n_1)\cdots\alpha(-n_r)F^m
+ u
\]\
where $u\in V_L^+(m)$ and $\ell_\alpha(u)<r$.
\end{lemma}

\subsection{Reduction I: even case}
\label{subsection:3.2}
In this subsection we prove by induction on $n$ that 
$V_L^+(n)\equiv 0\, \mod\, O(V_L^+)+M(1)^+$ for even integers $n$.

We need the following notation:
\begin{equation}\label{eqn:3.3}
\exp\left(
\sum_{n=1}^\infty\frac{x_n}{n}z^n
\right)
=\sum_{j=0}^\infty p_j(x_1,x_2,\dots)z^j = \sum_{j=0}^\infty p_j(x)z^j
\end{equation}
and $q_j(x) = p_j(x)+ p_j(-x)$. Then $p_j(x)$ are elementary 
Schur polynomials. 

The following lemma is easily derived from the definition of
vertex operators.

\begin{lemma}\label{lemma:3.4}
For $m,n\in\Z$,
\[
Y(e^{m\alpha},z)e^{n\alpha}
= \sum_{j=0}^\infty
p_j(m\alpha)\otimes e^{(m+n)\alpha} z^{2kmn+j}.
\]
where $p_j(\beta)=p_j(\b(-1),\b(-2),...).$ 
\end{lemma}

\begin{lemma}\label{lemma:3.5}
For any $m\in \Z_{>0}$,
\[
e^{2m\alpha}+e^{-2m\alpha}\equiv 0\quad\mod\quad O(V_L^+)+M(1)^+.
\]
\end{lemma}
\begin{proof}
Using Lemma \ref{lemma:3.4} we see that
\begin{align*}
&Y(E^m,z)E^m\\
&= Y(e^{m\alpha},z)e^{m\alpha}+Y(e^{-m\alpha},z)e^{-m\alpha}+
Y(e^{m\alpha},z)e^{-m\alpha}+Y(e^{-m\alpha},z)e^{m\alpha}\\
&=\sum_{j=0}^\infty
\left\{
p_j(m\alpha)\otimes e^{2m\alpha}+p_j(-m\alpha)\otimes e^{-2m\alpha}
\right\} z^{2km^2+j}
+ \sum_{j=0}^\infty q_j(m\alpha)z^{-2km^2+j}.
\end{align*}
Hence,
\[
\Res_{z}\frac{(1+z)^{km^2}}{z^{2km^2+1}}
Y(E^m,z)E^m
= E^{2m} + u
\]
where $u$ is a linear combination of $q_j(m\alpha)$,
in particular, $u\in M(1)^+$. Since $\wt (E^m)=km^2$ and
$2km^2+1\geq 2$ by Proposition \ref{proposition:2.4} (i),
$\Res_{z}\frac{(1+z)^{km^2}}{z^{2km^2+1}}
Y(E^m,z)E^m$ lies in $O(V_L^+).$ The proof is complete.
\end{proof}

\begin{lemma}\label{lemma:3.6}
Let $m\in\Z_{>0}$ be even and $n\in\Z_{>0}$. Then
\[
\alpha(-n)(e^{m\alpha}-e^{-m\alpha})\equiv 0\quad
 \mod\quad O(V_L^+)+M(1)^+.
\]
\end{lemma}
\begin{proof}
We prove the result by induction on $n$. Note that $L(-1)\sim -L(0)$.
Then
\[
L(-1)E^m \equiv -L(0)E^m \quad \mod \, O(V_L^+)
\]
or,
\[
m\alpha(-1)F^m
\equiv -km^2E^m
\quad \mod \, O(V_L^+).
\]
Using Lemma \ref{lemma:3.5} shows that 
$\alpha(-1)F^m\equiv 0\,\mod\, O(V_L^+)+M(1)^+.$ So the case $n=1$ is done.

Assume that the lemma is true for all integers less than $n$. 
Then by induction hypothesis,
\[
\alpha(-n)(e^{m\alpha}-e^{-m\alpha})\in M(1)^+\quad \mod\, O(V_L^+).
\]
Then again use $L(-1)\sim -L(0)$ to get
\begin{align}\label{eqn:3.4}
-(&n+km^2)\alpha(-n)F^m\notag\\
&=-L(0)\alpha(-n)F^m\notag\\
&\sim L(-1)(\alpha(-n)F^m)\notag\\
&=n\alpha(-n-1)F^m
+m\alpha(-n)\alpha(-1)E^m
\in M(1)^+\quad \mod\, O(V_L^+).
\end{align}

Let $v = \alpha(-n)\alpha(-1)\1$. Then by Lemma \ref{lemma:3.1} and 
induction hypothesis, we have
\begin{align*}
v*E^m
\equiv&
\alpha(-n)\alpha(-1)E^m\\
&+ 2mk(n+(-1)^{n+1})\alpha(-n-1)F^m\quad \mod\quad O(V_L^+)+M(1)^+.
\end{align*}
On the other hand, $E^m\in M(1)^+\,(\mod\, O(V_L^+)).$
Hence
\[
v*(e^{m\alpha}+e^{-m\alpha})\in M(1)^+\quad\mod\, O(V_L^+).
\]
and 
\[
\alpha(-n)\alpha(-1)E^m
\equiv -2mk(n+(-1)^{n+1})\alpha(-n-1)F^m\quad \mod\quad O(V_L^+)+M(1)^+.
\]
Finally, substituting this into (\ref{eqn:3.4}), we reach to
\[
\left\{
n-2km^2(n+(-1)^{n+1})
\right\}
\alpha(-n-1)F^m\equiv 0
\quad \mod\quad O(V_L^+)+M(1)^+.
\]
Since for $m\geq 2$ 
\[
n-2km^2(n+(-1)^{n+1})\neq 0,
\]
we see
\[
\alpha(-n-1)F^m\equiv 0
\quad \mod\quad O(V_L^+)+M(1)^+.
\]
\end{proof}

The main result in this subsection is the following:

\begin{lemma}\label{lemma:3.7} For any positive even integer $m,$ 
\[
V_L^+(m)\equiv 0\quad \mod\quad O(V_L^+)+M(1)^+.
\]
\end{lemma}
\begin{proof}
Let 
\[
u = \alpha(-n_1)\cdots \alpha(-n_r)(e^{m\alpha}+(-1)^re^{-m\alpha}).
\]
Set 
\[
v =
\begin{cases}
\alpha(-n_1)\cdots \alpha(-n_r)\1&\quad \text{if $r$ is even},\\
\alpha(-n_1)\cdots \alpha(-n_{r-1})\1&\quad \text{if $r$ is odd}
\end{cases}
\]
and 
\[
w =
\begin{cases}
e^{m\alpha}+e^{-m\alpha}&\quad \text{if $r$ is even},\\
\alpha(-n_r)(e^{m\alpha}-e^{-m\alpha})&\quad \text{if $r$ is odd.}
\end{cases}
\]
Then by Lemma \ref{lemma:3.2} and Lemma \ref{lemma:3.3}, we see that
$v*w = u + u'$ where $u'\in V_L^+(m)$ and $\ell_\alpha(u')<\ell_\alpha(u)= r$.
{}From Lemma \ref{lemma:3.5} and Lemma \ref{lemma:3.6},
$v*w \equiv 0\,\mod\, O(V_L^+)+M(1)^+,$ that is,
$u+u' \equiv 0\,\mod\, O(V_L^+)+M(1)^+$. An induction on $r$ shows that
$u\in M(1)^+\,\mod\, O(V_L^+)$.
\end{proof}

\subsection{Reduction II: odd case}
\label{subsection:3.3}

In this subsection  we prove that
\[
V_L^+(m)\subset M(1)^+\otimes (e^\alpha+e^{-\alpha})+O(V_L^+)
\]
for all odd integer $m$. Recall from the previous subsections that
\[
E=e^\alpha+e^{-\alpha},\quad F= e^\alpha-e^{-\alpha}.
\]

\begin{lemma}\label{lemma:3.8}
For any $s \in M(1)^+$ and $n\in\Z_{>0}$, 
there exists $t\in M(1)^+$ such that
\[
\alpha(-n)s\otimes F\equiv
t\otimes E\quad \mod\quad O(V_L^+).
\]
\end{lemma}
\begin{proof}
We  prove the lemma by induction on $n$. 
Applying $L(-1)$ to $s\otimes E$ and using the relation $L(-1)\sim L(0)$
yields
\begin{align*}
\alpha(-1)s\otimes F=&L(-1)(s\otimes E)-(L(-1)s)\otimes E\\
\equiv& -L(0)(s\otimes E)-(L(-1)s)\otimes E
\quad\mod\quad O(V_L^+).
\end{align*}
Since $M(1)^+\otimes (e^\alpha+e^{-\alpha})$ is invariant under $L(0)$
and $L(-1)s\in M(1)^+$ we see immediately that
\[
\alpha(-1)s\otimes F\in M(1)^+\otimes(e^\alpha+e^{-\alpha})
\quad \mod\quad O(V_L^+).
\]

Let us assume that the lemma holds for $n>0$. 
Again applying $L(-1)$ to $\alpha(-n)s\otimes F$ gives  
\begin{align*}
L(-1)(\alpha(-n)s&\otimes F)\\
&= (L(-1)s)\alpha(-n)\otimes F+
ns\alpha(-n-1)\otimes F
+s\alpha(-n)\alpha(-1)\otimes E.
\end{align*}
Thus
\begin{align*}
&ns\alpha(-n-1)\otimes F\\
&=L(-1)(s\alpha(-n)\otimes F)
- (L(-1)s)\alpha(-n)\otimes F
-s\alpha(-n)\alpha(-1)\otimes E\\
&\equiv -L(0)(s\alpha(-n)\otimes F)
-(L(-1)s)\alpha(-n)\otimes F
-s\alpha(-n)\alpha(-1)\otimes E
\quad \mod\quad O(V_L^+).
\end{align*}
By induction hypothesis both $L(0)(s\alpha(-n)\otimes F)$
and $(L(-1)s)\alpha(-n)\otimes F$ lie in $M(1)^+\otimes E$ modulo
$O(V_L^+).$ As a result we have
\[
\alpha(-n-1)s\otimes F
\in 
M(1)^+\otimes(e^\alpha+e^{-\alpha})
\quad\mod\quad O(V_L^+).
\]
\end{proof}

\begin{remark}\label{remark:3.9}
From the proof of Lemma \ref{lemma:3.8}, it is clear that for any $0\ne \gamma\in L$
\[
M(1)^-\otimes(e^\gamma-e^{-\gamma})\subset M(1)^+\otimes(e^\gamma+e^{-\gamma})
\quad\mod\quad O(V_L^+).
\]
\end{remark}

The main result in this subsection is:
\begin{lemma}\label{lemma:3.10} For any odd positive integer $m$, $V_L^+(m)\subset V_L^+(1)+O(V_L^+)$.
\end{lemma}
\begin{proof}
We prove this by induction on $m$. By Lemma \ref{lemma:3.8} 
the lemma holds for $m=1$. Suppose the assertion is true for 
$m-2\,(m\geq 3)$. A straightforward computation using Lemma \ref{lemma:3.4}
gives
\begin{align*}
&Y(E,z)E^{m-1}\\
&=Y(e^\alpha,z)e^{(m-1)\alpha}+Y(e^{-\alpha},z)e^{-(m-1)\alpha}+
Y(e^{-\alpha},z)e^{(m-1)\alpha}+Y(e^{\alpha},z)e^{-(m-1)\alpha}\\
&=\sum_{j=0}^\infty p_j(\alpha)\otimes e^{m\alpha}z^{2k(m-1)+j}
+\sum_{j=0}^\infty p_j(-\alpha)\otimes e^{-m\alpha}z^{2k(m-1)+j}\\
&\quad +\sum_{j=0}^\infty p_j(-\alpha)\otimes e^{(m-2)\alpha}z^{-2k(m-1)+j}
+\sum_{j=0}^\infty p_j(\alpha)\otimes e^{-(m-2)\alpha}z^{-2k(m-1)+j}\\
&= \sum_{j=0}^\infty\left\{
p_j(\alpha)\otimes e^{m\alpha}+p_j(-\alpha)\otimes e^{-m\alpha}
\right\}z^{2k(m-1)+j}\\
&+\sum_{j=0}^\infty \left\{p_j(-\alpha)\otimes e^{(m-2)\alpha} +p_j(\alpha)\otimes e^{-(m-2)\alpha}\right\}z^{-2k(m-1)+j}.
\end{align*}
So
\[
\Res_{z}\frac{(1+z)^k}{z^{2k(m-1)+1}}Y(E,z)
E^{m-1}
=e^{m\alpha}+e^{-m\alpha}+u
\]
where $u\in V_L^+(m-2)$. Since $\wt (E)=k$ and $2k(m-1)+1\geq 2$
we see from Proposition \ref{proposition:2.4} that 
$\Res_{z}\frac{(1+z)^k}{z^{2k(m-1)+1}}Y(E,z)
E^{m-1}\in O(V_L^+).$ From the induction hypothesis we know that
$u$ lies $V_L^+(1)$ modulo $O(V_L^+).$ Thus
\[
e^{m\alpha}+e^{-m\alpha}\in V_L^+(1)\quad\mod\quad O(V_L^+).
\]
By the same argument given in the proof of Lemma \ref{lemma:3.6} 
we prove that
\[
\alpha(-n)\otimes F^m\in V_L^+(1)\quad\mod\quad O(V_L^+).
\]

Next, we want to show by induction on $r$ that
\[
u = \alpha(-n_1)\alpha(-n_2)\cdots \alpha(-n_r)(e^{m\alpha}
+(-1)^r e^{-m\alpha})
\in V_L^+(1)+O(V_L^+).
\]
Set 
\[
v =
\begin{cases}
\alpha(-n_1)\cdots \alpha(-n_r)\1&\quad \text{if $r$ is even},\\
\alpha(-n_1)\cdots \alpha(-n_{r-1})\1&\quad \text{if $r$ is odd}
\end{cases}
\]
and
\[
w =
\begin{cases}
e^{m\alpha}+e^{-m\alpha}&\quad \text{if $r$ is even},\\
\alpha(-n_r)(e^{m\alpha}-e^{-m\alpha})&\quad \text{if $r$ is odd.}
\end{cases}
\]
Then from Lemmas \ref{lemma:3.2} and \ref{lemma:3.3} we see that
\[
v*w = \alpha(-n_1)\alpha(-n_2)\cdots\alpha(-n_r)(e^{m\alpha}
+(-1)^r e^{-m\alpha})+u'
\]
where $\ell_{\alpha}(u')<r$ and $u'\in V_L^+(m)$.
By induction hypothesis  $w,u'\in V_L^+(1)+O(V_L^+).$ Since $v\in M(1)^+$
we have $v*w \in V_L^+(1)+O(V_L^+)$ and 
then $u\in V_L^+(1)+O(V_L^+).$
\end{proof}

\section{Generators of  $A(V_L^+)$}
\setcounter{equation}{0}
\label{section:4}

We have already proved in Section 3  that $A(V_L^+)=M(1)^++V_L^+(1)+O(V_L^+).$ 
The main result of this section is 
that $A(V_L^+)$ is generated by $\o+O(V_L^+),J+O(V_L^+)$ and $E+O(V_L^+).$ 
Since $M(1)^++O(V_L^+)$ is generated by $\o+O(V_L^+)$ and $J+O(V_L^+)$
[DN] we establish that $V_L^+(1)+O(V_L^+)$ is generated by $\o+O(V_L^+)$
and $E+O(V_L^+).$ Since the structure of $V_L^+(1)$ as a Virasoro module varies 
according to whether $k$ is a perfect square or not, we deal with
these cases separately.  
The case that $k$ is a perfect is more complicated.
Nevertheless, the ideas and
the techniques developed in [DN] still work in the present situation. 

\subsection{A spanning set of $V_L^+(1)+O(V_L^+)$ I: $k$ is not a perfect 
square}
\label{subsection:4.1}

In this section we assume that $k$ is not a perfect square.
In this case $V_L^+(1)$ is an irreducible Virasoro module 
which is isomorphic to $L(1,k)$
with a highest weight vector $e^\alpha+e^{-\alpha}$.
 For short, we set 
\[
v^{*s} = \overset{s}{\overbrace{v*\cdots*v}}
\]
for $v\in V_L^+$. Recall that $[v]=v+O(V_L^+)$ for $v\in V_L^+$,
we will use a similar notation $[v]^{*s}$. Then it is easy to see that
$[v^{*s}] = [v]^{*s}$.

\begin{lemma}\label{lemma:4.1}
Suppose that $k$ is not a perfect square. 
Then $V_L^+(1)+O(V_L^+)$ is spanned by
\[
[\mathcal{S}_{\omega,\,E}]=
\left\{
[\omega^{*s}*E]\,|\,s\geq 0
\right\}.
\]
\end{lemma}
\begin{proof} In this case $V_L^+(1)$ is spanned by the vectors 
\[
v = L(-n_1)L(-n_2)\cdots L(-n_r)E,\quad
n_1\geq n_2\geq\dots\geq n_r\geq 1.
\]
So it is enough to show that $[v]$ is spanned
by $\mathcal{S}_{\omega,\,E}$. Using Proposition \ref{proposition:2.4}
(iii), (iv) and the relation
\[
L(0)L(-n_1)\cdots L(-n_r)E= (n_1+\cdots+n_r+k)L(-n_1)\cdots L(-n_r)E
\]
one can easily show that
$[v] = [P(\omega)*E]$ with some polynomial $P(x)$.
\end{proof}

\subsection{A spanning set of $V_L^+(1)+O(V_L^+)$ II: $k$ is a perfect 
square}
\label{subsction:4.2}

In this subsection  we consider the case that $k$ is a perfect square.
Since $V_L^+(1)$ in an irreducible $M(1)^+$ module and $M(1)^+$ is
generated by $\o$ and $J$ one can see that $V_L^+(1)$ 
is spanned by
\[
\left\{
u^1_{m_1}\cdots u^k_{m_k}E\,|\,u^i=\omega, J,\, m_i\in \Z
\right\}
\]
which are not necessarily linearly independent. 
We say that an expression $u^1_{m_1}\cdots u^k_{m_k}E$ has length $t$ 
with respect to $J$,  
which we write $\ell_J(u^1_{m_1}\cdots u^k_{m_k}E) = t,$
if $\{i|u^i=J\}$ has cardinality $t.$ Note that $\omega_i=L(i-1).$ 
An induction on $\ell_J(u^1_{m_1}\cdots u^k_{m_k}E)$ 
using Lemma \ref{lemma:3.1new} (1)  shows that 
$u^1_{m_1}\cdots u^k_{m_k}E$ is a linear combination of vectors 
of type
\begin{align*}
&\left\{
L(m_1)L(m_2)\cdots L(m_s)J_{n_1}J_{n_2}\cdots J_{n_t}E\,|\,
m_a,n_b\in\Z
\right\}.
\end{align*}

Using commutation relation in Lemma \ref{lemma:3.1new} and
the fact that $E$ is a singular vector we can prove the following lemma.

\begin{lemma}\label{lemma:4.2new}
Let $W$ be a subspace of $V_L^+$ spanned by $J_{n_1}\cdots J_{n_t}E$
with $n_i\in \Z$. Then $W$ is invariant under the action of 
$L(m),\,m\geq 0$.
\end{lemma}

\begin{lemma}\label{lemma:4.2}
$V_L^+(1)$ is spanned by
\[
L(-m_1)\cdots L(-m_s)J_{-n_1}\cdots J_{-n_t}E
\]
where $m_1\geq m_2\geq \dots\geq m_s\geq 1$ and
$n_1\geq n_2\geq\dots\geq n_t\geq 1$.
\end{lemma}
\begin{proof} 
 We have already known that $V_L^+(1)$ is spanned by
\[
L(-m_1)\cdots L(-m_s)J_{-n_1}\cdots J_{-n_t}E
\]
where $m_a,n_b\in\Z$. Using the PBW theorem for the Virasoro algebra
we can assume that $m_1\geq \cdots \geq m_s.$ By 
Lemma \ref{lemma:4.2new} we can further assume that
$m_1\geq m_2\geq \dots\geq m_s\geq 1$.  
We proceed by induction on $\ell_J(v)$ 
that $v=L(-m_1)\cdots L(-m_s)J_{-n_1}\cdots J_{-n_t}E$
can be spanned by the indicated vectors in the proposition.

If the length is $0$, it is clear.
Suppose that it is true for all monomials $v$ such that 
$\ell_J(v)<t$. By lemma \ref{lemma:3.2new} and
the induction hypothesis we can assume $n_t\geq 1.$
If $n_1\geq \cdots \geq n_t$ we are done. Otherwise 
there exists $n_a$ such that $n_{a+1}\geq \cdots \geq n_t$ but
$n_a< n_{a+1}.$ There are two cases $n_a\leq 0$ and $n_a >0$ 
which are dealt with separately.
If $n_a\leq 0$, then  
\begin{align*}
L(-m_1)\cdots&L(-m_s)J_{-n_1}\cdots J_{-n_t}E\\
&= \sum_{j=a+1}^t L(-m_1)\cdots L(-m_s)J_{-n_1}\cdots
\overset{\,\,\,\vee}{J}_{-n_a}
\cdots [J_{-n_a},J_{-n_j}]\cdots J_{-n_t}E\\
&\quad\quad 
+\sum_{j=a+1}^t L(-m_1)\cdots L(-m_s)J_{-n_1}\cdots
\overset{\,\,\,\vee}{J}_{-n_a}\cdots J_{-n_t}J_{-n_a}E
\end{align*}
where $\overset{\,\,\,\vee}{J}_{-n_a}$ means that we omit the term 
$J_{-n_a}$.
However by Lemma \ref{lemma:3.1new} (2), $[J_{-n_a}, J_{-n_j}]$ are 
linear combinations of operators of type
$$L(p_1)\cdots L(p_{s'}),\quad L(q_1)\cdots L(q_{t'})J_r.$$
By substituting these into the above and using commutation relation
in Lemma \ref{lemma:3.1new} (1) again, the first term of the right hand side 
is a linear combination of monomials whose lengths with respect to $J$ are 
less than or equal to $t-1$. Further by lemma \ref{lemma:3.2new},
the second term is also a linear combination of such monomials.
Thus by induction hypothesis,
this is expressed as linear combinations of expected monomials. 

If $n_a>0$ then either $n_a<n_t$ or there exists $b$ with $t>b>a$ so that
$n_b>n_a \geq n_{b+1}.$ Then we have either  
\begin{align*}
L(-m_1)\cdots& L(-m_s)J_{-n_1}\cdots J_{-n_t}E\\
&= \sum_{j=a+1}^t L(-m_1)\cdots L(-m_s)J_{-n_1}\cdots
\overset{\,\,\,\vee}{J}_{-n_a}
\cdots [J_{-n_a},J_{-n_j}]\cdots J_{-n_t}E\\
&\quad\quad\quad+ L(-m_1)\cdots L(-m_s)J_{-n_1}\cdots
\overset{\,\,\,\vee}{J}_{-n_a}
\cdots J_{-n_t}J_{-n_a}E
\end{align*}
or
\begin{align*}
L(-m_1)\cdots& L(-m_s)J_{-n_1}\cdots J_{-n_t}E\\
&= \sum_{j=a+1}^b L(-m_1)\cdots L(-m_s)J_{-n_1}\cdots
\overset{\,\,\,\vee}{J}_{-n_a}
\cdots [J_{-n_a},J_{-n_j}]\cdots J_{-n_t}E\\
&\quad\quad+ L(-m_1)\cdots L(-m_s)J_{-n_1}\cdots
\overset{\,\,\,\vee}{J}_{-n_a}
\cdots J_{-n_b}J_{-n_a}J_{-n_{b+1}}\cdots J_{-n_t}E
\end{align*}
{}From the discussion of case $n_a\leq 0$ it is enough to show
either 
\[
 L(-m_1)\cdots L(-m_s)J_{-n_1}\cdots\overset{\,\,\,\vee}{J}_{-n_a}
\cdots J_{-n_t}J_{-n_a}E
\]
or 
\[
 L(-m_1)\cdots L(-m_s)J_{-n_1}\cdots\overset{\,\,\,\vee}{J}_{-n_a}
\cdots J_{-n_b}J_{-n_a}J_{-n_{b+1}}\cdots J_{-n_t}E
\]
can be expressed as linear combinations of desired vectors. But this
follows from an induction on $a.$ 
\end{proof}

\begin{lemma}\label{lemma:4.3} Assume that $k\ne 1.$
$V_L^+(1)+O(V_L^+)$ is spanned by
\[
[\mathcal{S}_{\omega,\,E}] =\left\{
[\omega^{*s}*E]\,|\, s\in \Z_{\geq 0}
\right\}.
\]
\end{lemma}
\begin{proof} The case that $k$ is not a perfect square was treated in Lemma \ref{lemma:4.1} already. So we can assume that $k$ is a perfect square. By lemma \ref{lemma:4.2},
it is enough to show that
any 
\[
[v]= [L(-m_1)L(-m_2)\cdots L(-m_s)
J_{-n_1}\cdots J_{-n_t}E]
\]
where $m_1\geq m_2\geq \dots\geq m_s\geq 1$ and 
$n_1\geq n_2\geq \dots\geq n_t\geq 1$ is spanned by $[\mathcal{S}_{\omega,\,E}]$.
We prove this by induction on $\ell_J(v)$. If the length is $0$, 
the proof of Lemma \ref{lemma:4.1} gives the result. 

Let $t>0$ and assume that
the statement is true for all $v$ with $\ell_J(v)<t$. We will prove that
$[v]$ is spanned by $[\mathcal{S}_{\omega,\,E}]$ by induction on weight of $v$.
Clearly, the smallest weight is $4t+k$ and $v = J_{-1}\cdots J_{-1}E$. Then
\[
J*\cdots *J*E
- v
=\sum_{
\begin{Sb}
n_i\in\{-1,0,1,2,3\},\\
(n_i)\neq (-1,-1,\dots,-1)
\end{Sb}
}
a_{n_1\dots n_t}J_{n_1}\cdots J_{n_t}E.
\]
Since each term appeared in the right hand side involves  $J_{n_i}$ for
some nonnegative integer $n_i$, by using Lemma \ref{lemma:3.2new},
it length is strictly less than $t$. 
Thus by induction hypothesis, the image of right hand side in $A(V_L^+)$
is spanned by
$[\mathcal{S}_{\omega,\,E}].$ So we can assume that $v=J*\cdots *J*E.$
Note that
\[
J*E= \sum_{j=0}^4\binom{4}{j}J_{j-1}E
\]
and $\wt(J_{j-1}E) = 4+k-j\leq 4+k$ if $j\geq 0$. Then from the
decomposition of $V_L^+(1)$ (see equation (\ref{eqn:3.2new})) we see that
$J*E$ is a vector in the irreducible module for the Virasoro algebra
generated by the highest weight $E.$ The proof of Lemma \ref{lemma:4.1}
shows that $[J*E]$ is a linear combination of elements of 
$[\mathcal{S}_{\omega,\,E}]$. Then using the fact that
$\omega$ is a central element proves that $v$ is spanned 
by $[\mathcal{S}_{\omega,\,E}]$.

Now consider a general vector
\[
v = L(-m_1)L(-m_2)\cdots L(-m_s)
J_{-n_1}\cdots J_{-n_t}E.
\]
Suppose $m_1>2$. Then by using Proposition \ref{proposition:2.4} (3)
we have
\[
v\sim (-1)^{m_1}\left\{
(m_1-1)(L(-2)+L(-1))+L(0)
\right\}L(-m_2)\cdots L(-m_s)J_{-n_1}\cdots J_{-n_t}E
\]
which is a sum of three homogeneous vectors of weight strictly less than 
$\wt(v)$. Then by induction hypothesis, $[v]$ is spanned by 
$[\mathcal{S}_{\omega,\,E}]$. Thus we can assume that $m_1\leq 2.$ 
We can further assume by using the relation $L(-1)\sim -L(0)$ that
$m_1=m_2=\dots=m_s=2$. Namely,
\[
v =\overset{s}{\overbrace{L(-2)\cdots L(-2)}}
J_{-n_1}J_{-n_2}\cdots J_{-n_t}E.
\]
Then
\[
v=
\omega^{*s}*(J_{-n_1}\1)*
(J_{-n_2}\cdots J_{-n_t}\1)
*E+u\]
where the weights of homogeneous components of $u$ are less than $\wt(v)$
 and the length of each homogeneous component of $u$ with respect to $J$ 
is less than or equal to $t.$ 
Then again by using induction hypothesis $[u]$ is spanned by
$[\mathcal{S}_{\omega,\,E}]$. It reduces to the case that
\[
v =\omega^{*s}*(J_{-n_1}\1)*
(J_{-n_2}\cdots J_{-n_t}\1)
*E.
\]
Note that
\begin{align*}
v&\equiv (J_{-n_1}\1)*\omega^{*s}*(J_{-n_2}\cdots J_{-n_t}\1)*E\\
&=(J_{-n_1}\1)*\omega^{*s}*
(J_{-n_2}\cdots J_{-n_t}E)+w
\end{align*}
where $w$ is a vector spanned by $[\mathcal{S}_{\omega,\,E}]$.
By induction hypothesis of the length with respect to $J,$ we see
that $J_{-n_2}\cdots J_{-n_t}E$ is spanned by $\mathcal{S}_{\omega,\,E}$
modulo $O(V_L^+).$ So we can assume that 
$$v=\o^{*p}*(J_{-n_1}\1)*E$$
for some $p\geq 0.$ Again by induction hypothesis on the length
of $v$ with respect to $J$ we conclude that such $[v]$ is spanned by
$[\mathcal{S}_{\omega,\,E}].$ This establishes the lemma.
\end{proof}

Let us  summarize the main results in this section.

\begin{proposition}\label{proposition:4.4} Assume that $k\ne 1.$ Then the 
Zhu's algebra $A(V_L^+)$ is spanned by
$$\{[\o^{*s}*J^{*t}], [\o^{*s}*E]|s,t\geq 0\}$$
\end{proposition}

\section{The structure of $A(V_L^+)$}
\label{section:5}

It is proved in Section 4 that the algebra $A(V_L^+)$ is 
generated by $[\o],$ $[J]$ and $[E]$ if $k\ne 1.$ 
In this section we determine the algebra structure of $A(V_L^+)$
which is  a commutative semisimple algebra of dimension $k+7$ if
$k\ne 1.$  
This is achieved by studying the relations 
among $[\omega], [J]$ and $[E].$ We have already known two relations
between $[\o]$ and $[J]$ from [DN]. Using the known irreducible
modules of $A(V_L^+),$ we obtain more relations. The classification
of irreducible modules for $V_L^+$  follows immediately from
the dimension of $A(V_L^+)$ as $A(V_L^+)$ has $k+7$ known
irreducible modules. 

In Subsection 5.1 we list all known irreducible modules of $V_L^+$
and give the scalars of $\o,J,E$ on the top levels of these
modules. In Subsection 5.2 we find two relations among
$[\o],$ $[J]$ and $[E].$ Subsection 5.3 is the core of this
paper where we determine a basis of $A(V_L^+).$  
Subsection 5.4 is easy but important.
In this subsection we classify the irreducible modules for $V_L^+.$

We assume that $k\ne 1$ in the first three subsections.
\subsection{List of known irreducible modules}
\label{subsection:5.1}

Here we give the list of known irreducible $V_L^+$-modules and the action
of $\omega, E$ and $J$ on the top levels of them.
As we mentioned before, we have the following irreducible
$V_L^+$-modules
\begin{equation}
\begin{split}
&V_L^+,\quad V_L^-,\quad V_{L+\frac{r}{2k}\alpha}\,(r= 1,2,\dots,k-1),\\
&V_{L+\frac{\alpha}{2}}^+,\quad V_{L+\frac{\alpha}{2}}^-,
\quad V_L^{T_1,+}, \quad V_L^{T_1,-},\quad V_L^{T_2,+},\quad V_L^{T_2,-}.
\end{split}
\end{equation}
Note that the top level of these modules are 1-dimensional and
$\o$, $J$ and $E$ act as scalars. 
The following table gives the scalars which follows from the 
construction of these modules (cf. [DN]). 
\bigskip

\begin{center}

\begin{tabular}{|c|c|c|c|c|c|}\hline
&$V_L^+$&$V_L^-$&$V_{L+\frac{r}{2k}\alpha}\,(1\leq r\leq k-1)$&$
V_{L+\frac{\alpha}{2}}^+$&$V_{L+\frac{\alpha}{2}}^-$\\ \hline
$\omega$&$0$&1&$r^2/4k$&$k/4$&$k/4$\\ \hline
$E$&0&0&0&1&$-1$\\ \hline
$J$&0&$-6$&$c^4-c^2/2, c^2 = r^2/2k$&$k^4/4-k/4$&$k^4/4-k/4$\\ \hline
\end{tabular}

\bigskip
\begin{tabular}{|c|c|c|c|c|}\hline
&$V_L^{T_1,+}$&$ V_L^{T_1,-}$&$V_L^{T_2,+}$&$V_L^{T_2,-}$\\ \hline
$\omega$&$1/16$&$9/16$&$1/16$&$9/16$\\ \hline
$E$&$2^{-2k+1}$&$-2^{-2k+1}(4k-1)$&
$-2^{-2k+1}$&$2^{-2k+1}(4k-1)$\\ \hline
$J$&3/128&$-45/12$8&
3/128&$-45/128$\\ \hline
\end{tabular}

\end{center}

\subsection{The relations among $\omega,$ $E$ and $J$}
In this subsection  we first  prove the relation
\begin{equation}
([\omega]-k/4)*([\omega]-1/16)*([\omega]-9/16)*[E] = 0.
\end{equation}
Note that $V_L^+(1)=\oplus_{n\geq k}V_L^+(1,n)$ is $\Z$-graded 
where $V_L^+(1,n)$ is the weight $n$ subspace of $V_L^+(1).$ 
Then  
$$(\omega-k/4)*(\omega-1/16)*(\omega-9/16)*E\in  \oplus_{0\leq n\leq k+6}V_L^+(1,n).$$

It is easy to see that $V_L^+(1,k+6)$ has  the following basis:
\[
\begin{array}{ll}
g_1= \alpha(-6)F,&\quad g_2 = \alpha(-5)\alpha(-1)E\\
g_3=\alpha(-4)\alpha(-2)E,&\quad g_4=\alpha(-4)\alpha(-1)^2F,\\
g_5 = \alpha(-3)^2E,&\quad 
g_6 = \alpha(-3)\alpha(-2)\alpha(-1)F,\\
g_7 = \alpha(-3)\alpha(-1)^3E,&\quad
g_8 = \alpha(-2)^3F,\\
g_9 = \alpha(-2)^2\alpha(-1)^2E,&\quad
g_{10} = \alpha(-2)\alpha(-1)^4F,\\
g_{11} = \alpha(-1)^6E.&
\end{array}
\]
In particular, $\dim V_L^+(k,6) = 11$. 
Similarly $V_L^+(1,k+5)$ has the following basis  
\[
\begin{array}{ll}
f_1=\alpha(-5)F,&\quad f_2 = \alpha(-4)\alpha(-1)E,\\
f_3=\alpha(-3)\alpha(-2)E,&\quad
f_4 = \alpha(-3)\alpha(-1)^2F,\\
f_5=\alpha(-2)^2\alpha(-1)F,&\quad
f_6= \alpha(-2)\alpha(-1)^3E,\\
f_7=\alpha(-1)^5F&
\end{array}
\]
and dimension 7. We also need  the following  basis of 
$V_L^+(1,k+3)$:
\[
h_1= \alpha(-3)F,\quad 
h_2 = \alpha(-2)\alpha(-1)E,\quad
h_3=\alpha(-1)^3 F.
\]
Clearly, $\dim\, V_L^+(k,3)= 3$.

Set $v = \alpha(-1)^4_{-3}E$ where $\a(-1)^4_{-3}$ is
the component operator of $Y(\a(-1)^4,z)=\sum_{n\in\Z}\a(-1)^4_nz^{-n-1}.$ 

\begin{lemma}\label{lemma:5.1}
The vectors
\[
L(-1)(f_i)\quad(i=1,\dots,7),
\quad
L(-3)(h_j)\quad (j= 1,2,3),
\quad
v
\]
form a basis of $V_L^+(1,k+6)$.
\end{lemma}

\begin{proof} The main idea of the proof is to show that these vectors are 
linearly independent. This is done in the following table 
by giving explicit expressions
of these vectors in terms of $g_i$ for $i=1,...,11.$ In fact if we denote 
the matrix below by $A$ then $\det A=6144(1-k)k^2.$ Thus $A$ is 
non-singular if $k\ne 1.$ 
\end{proof} 

\begin{center}
\begin{tabular}{|c|c|c|c|c|c|c|c|c|c|c|c|}
\hline
&$g_1$&$g_2$&$g_3$&$g_4$&$g_5$&$g_6$&$g_7$&$g_8$&$g_9$&$g_{10}$&$g_{11}$\\
\hline
$L(-1)f_1$&5&1&0&0&0&0&0&0&0&0&0\\
\hline
$L(-1)f_2$&0&4&1&1&0&0&0&0&0&0&0\\
\hline
$L(-1)f_3$&0&0&3&0&2&1&0&0&0&0&0\\
\hline
$L(-1)f_4$&0&0&0&3&0&2&1&0&0&0&0\\
\hline
$L(-1)f_5$&0&0&0&0&0&4&0&1&1&0&0\\
\hline
$L(-1)f_6$&0&0&0&0&0&0&2&0&3&1&0\\
\hline
$L(-1)f_7$&0&0&0&0&0&0&0&0&0&5&1\\
\hline
$2kL(-3)h_1$&$6k$&0&0&0&$2k$&1&0&0&0&0&0\\
\hline
$2kL(-3)h_2$&0&$4k$&$2k$&0&0&$2k$&0&0&1&0&0\\
\hline
$2kL(-3)h_3$&0&0&0&$6k$&0&0&$2k$&0&0&1&0\\
\hline
$v$&$32k^3$&$48k^2$&$48k^2$&$24k$&$24k^2$&$48k$&4&$8k$&6&0&0\\
\hline

\end{tabular}
\end{center}

\begin{lemma}\label{lemma:5.2}
We have the relation
\[
(\omega-k/4)(\omega-1/16)(\omega-9/16)E = 0
\]
in $A(V_L^+)$ where we also use $v$ to denote its image in $A(V_L^+)$
for any $v\in V_L^+$ and the product is the $*$ operation.
\end{lemma}
\begin{proof}
We first note that all vectors of a basis in Lemma \ref{lemma:5.1} 
are congruent to the vectors of $V_L^+(1)$ of weight less than
or equal to $k+5.$ By Lemma \ref{lemma:4.3} there exists a monic polynomial 
of degree 3
such that
\begin{equation}\label{eqn:5.3}
f(\omega)E= 0.
\end{equation}
Let us apply both sides of (\ref{eqn:5.3}) to the top levels 
of $V_{L+\frac{\alpha}{2}}^+, V_L^{T_1,+}, V_L^{T_2,+}$ and note
that $E$ is nonzero on these top levels. We immediately have
\[
f(k/4) = f(1/16)=f(9/16) = 0.
\]
Since $f$ has degree $3$ we get
\[
f(x) = (x-k/4)(x-1/16)(x-9/16).
\]
\end{proof}

Next we study relations between $J$ and $E.$ 
{}From Lemma \ref{lemma:4.3}, it is clear that there exists a polynomial
$r(x)$ of degree 2 such that $J*E= r(\omega)E.$ We give the explicit
expression of the $r(x)$ in the following lemma.
\begin{lemma}
We have
\begin{equation}\label{eqn:5.4}
J*E =E*J= r(w)E
\end{equation}
in $A(V_L^+)$ where
\begin{equation}\label{eqn:5.5}
r(x) =
\frac{2(32k^2-8k-9)}{(4k-9)(4k-1)}x^2
+\frac{9+80k-104k^2}{2(4k-1)(4k-9)}x
+\frac{27k(k-1)}{8(4k-1)(4k-9)}.
\end{equation}
\end{lemma}
\begin{proof}
Set $r(x) = ax^2+bx+c$. We will evaluate $J*E=r(\o)E$ on the top levels
of  modules listed in Subsection \ref{subsection:5.1} on which 
$E\neq 0.$ Namely, we calculate the values of 
$\omega$ and $J$ on the top levels of modules $V_{L+\frac{\alpha}{2}}^+,
V_L^{T_1,+}$ and $V_L^{T_1,-}$. The we have
\begin{align*}
&k^2 a+4kb+16c= 4k^2-4k,\\
&a+16b+256c=6,\\
&81a+144b+256c= -90.
\end{align*}
By solving this linear system, we have the desired result. The same argument
also shows that $E*J= r(w)E.$  In particular, $J$ and $E$ are commutative.
\end{proof}

\subsection{A basis for $A(V_L^+)$}

So far, we have established the following relations;
\begin{align}
&J^2= p(\omega)+q(\omega)J,\tag{$B_1$}\\
&(\omega-1)(\omega-1/16)(\omega-9/16)(J+\omega-4\omega^2)= 0,\tag{$B_2$}\\
&JE= r(\omega)E,\tag{$L_1$}\\
&t(\omega)E= 0\tag{$L_2$}
\end{align}
where
\begin{align}
&p(x)=\frac{1816}{35}x^4-\frac{212}{5}x^3 +\frac{89}{10}x^2-\frac{27}{70}x,\\
& q(x)=-\frac{314}{35}x^2+\frac{89}{14}x-\frac{27}{70},\\
&r(x) =
\frac{2(32k^2-8k-9)}{(4k-9)(4k-1)}x^2
+\frac{9+80k-104k^2}{2(4k-1)(4k-9)}x
+\frac{27k(k-1)}{8(4k-1)(4k-9)},
\\
&t(x)= (x-k/4)(x-1/16)(x-9/16).
\end{align}
We remark that the relation ($B_1$) and ($B_2$) were found in [DN] in the
algebra $A(M(1)^+).$ Since $O(M(1)^+)\subset O(V_L^+)$ these two relations
are  also true in $A(V_L^+).$   

\begin{lemma}\label{lemma:5.4}
\[
E*E = \sum_{j=0}^k\binom{k}{j}q_{2k-j}(\alpha)\1.
\]
\end{lemma}
\begin{proof}
{}From the proof of Lemma \ref{lemma:3.5}, we have
\[
Y(E,z)E
= \sum_{j=0}^\infty
\left\{
p_j(\alpha)\otimes e^{2\alpha}+p_j(-\alpha)\otimes e^{-2\alpha}
\right\}
z^{2k+j}
+\sum_{j=0}^\infty q_j(\alpha)z^{-2k+j}.
\]
Therefore, we have
\[
E*E= \Res_{z}\left(
\frac{(1+z)^k}{z}Y(E,z)E
\right)
=\sum_{j=0}^k\binom{k}{j}E_{j-1}E
=\sum_{j=0}^k\binom{k}{j}q_{2k-j}(\alpha)\1.
\]
\end{proof}

\begin{lemma}\label{lemma:5.5}
There exist polynomials $a(x),\deg\, a = k$ and $s(x),\deg\,s\leq 2$
such that
\begin{equation}\label{eqn:5.10}
E^2= a(\omega)+s(\omega)(J+\omega-4\omega^2).
\end{equation}
Further
\[
a(x)= a_0x(x-\frac{1}{4k})(x-\frac{4}{4k})\cdots (x-\frac{(k-1)^2}{4k})
\]
where $a_0 = 2(4k)^k/(2k)!$.
\end{lemma}
\begin{proof}
Since $E*E\in M(1)^+$ and the highest weight of homogeneous component
is $2k$ we can write $E^2$ as a linear combination
of $\o^{*s}*J^{\*t}$ for $s,t\geq 0$  such that $2s+4t\leq 2k$ (see
[DN]). The existence of $a(x)$ and $s(x)$ follow from the relation
($B_1$)-($B_2$). Clearly, the degrees of $a(x)$ and $s(x)$ are less than
or equal to $k$ and $k-2,$ respectively. Using ($B_2$) we can assume
that the degree of $s(x)$ is less than or equal to 2.

Let us apply both hand sides of (\ref{eqn:5.10}) to the top level of modules
\[
V_L^+,V_{L+\frac{1}{2k}\alpha},\dots,V_{L+\frac{k-1}{2k}\alpha}.
\]
Since both $E$ and $J+\omega-4\omega^2$ act trivially
on these top levels, we see that $a(\omega)$ also acts trivially on 
top levels. This implies
\[
a(0)=a(\frac{1}{4k})=\dots=a(\frac{(k-1)^2}{4k})= 0.
\]
Since $\deg\,a \leq k$, we find 
\[
a(x)= a_0x(x-\frac{1}{4k})(x-\frac{4}{4k})\cdots (x-\frac{(k-1)^2}{4k})
\]
for some $a_0\in\C$. Note that $J+\omega-4\omega^2$ acts 
trivially and $E=1$ on the top level of 
$V_{L+\frac{\alpha}{2}}^+$, we have
$a(k/4) = 1$, 
which implies $a_0 = 2(4k)^k/(2k)!$.
\end{proof}

Set
\begin{equation}\label{ex3}
\varphi(x) =
(x-1)(x-\frac{1}{16})(x-\frac{9}{16})(x-\frac{k}{4})a(x).
\end{equation}

\begin{lemma}\label{lemma:5.6}
We have
\[
\varphi(\omega)= 0.
\]
\end{lemma}
\begin{proof}
From ($B_2$)
\[
(\omega-1)(\omega-\frac{1}{16})(\omega-\frac{9}{16})(J+\omega-4\omega^2) = 0
\]
we see from Lemma \ref{lemma:5.5} that
\[
(\omega-1)(\omega-\frac{1}{16})(\omega-\frac{9}{16})E^2
=(\omega-1)(\omega-\frac{1}{16})(\omega-\frac{9}{16})a(\omega).
\]
On the other hand, ($L_2$) tells us
\[
(\omega-\frac{k}{4})(\omega-\frac{1}{16})(\omega-\frac{9}{16})E^2=0
\]
and therefore $\varphi(\omega) = 0$.
\end{proof}

Since $JE= r(\omega)E$, we see
\begin{align*}
0&=(J-r(\omega))E^2\\
&=(J-r(\omega))a(\omega)+(J-r(\omega))(J+\omega-4\omega^2)s(\omega)
\end{align*}
and therefore
\begin{align*}
J^2s(\omega)+\{
a(\omega)+(\omega-4\omega^2-&r(\omega))s(\omega)\} J\\
&-r(\omega)a(\omega)-r(\omega)(\omega-4\omega^2)s(\omega)= 0.
\end{align*}
By using the relation $J^2 = p(\omega)+q(\omega)J$, this is reduced to
\begin{align}\label{ex1}
\{a(\omega)+(q(\omega)-r(\omega)+\omega-4\omega^2)s(\omega)\}&J
-r(\omega)a(\omega)\\
&+\left\{
p(\omega)-r(\omega)(\omega-4\omega^2)\right\}
s(\omega) = 0.\nonumber
\end{align}

For convenience we introduce   
\begin{equation}\label{ex2}
b(x)= a(x)+(q(x)-r(x)+x-4x^2)s(x).
\end{equation}
\begin{lemma}\label{lemma:5.8}
\[
b(1)= a(1)+\frac{27(-12+65k-33k^2)}{8(9-40k+16k^2)}s(1),
\quad b(\frac{1}{16}) = a(\frac{1}{16}),
\quad b(\frac{9}{16})=a(\frac{9}{16}).
\]
\end{lemma}
\begin{proof} A straightforward calculation shows that
\[
q(x)-r(x)+x-4x^2 = \frac{9(-12+65k-33k^2)(16x-1)(16x-9)}{280(4k-1)(4k-9)}.
\]
The lemma follows.
\end{proof}

\begin{lemma}\label{lemma:5.10} (1) If $k$ is not a perfect square, then 
$$b(1)\ne 0, b(\frac{1}{16})\ne 0,   b(\frac{9}{16})\ne 0.$$

\noindent
(2) If $k=4m^2$ for some positive integer $m$ then  
$$b(1)=b(\frac{1}{16})=b(\frac{9}{16})=0.$$

\noindent
(3) If $k=(2m+1)^2$ for some positive integer $m$ then 
$$b(1)=0, b(\frac{1}{16})\ne 0,   b(\frac{9}{16})\ne 0.$$
\end{lemma}

\begin{proof} (1) Recall that $a(x)=a_0\prod_{r=0}^{k-1}(x-\frac{r^2}{4k}).$
If $b(1/16)=a(1/16)=0,$ then there 
exists $i\,(1\leq i\leq k-1)$ such that $\frac{1}{16} = \frac{i^2}{4k},$
namely, $k= 4i^2$. This is a contradiction. 
 By the exactly same reason, 
we know $b(9/16)=a(9/16)\neq 0$. 
It remains to show $b(1)\neq 0$. Let us evaluate the relation
\[
E^2 = a(\omega)+s(\omega)(J+\omega-4\omega^2)
\]
on the top level of the module $V_L^-.$ Then we have $a(1) = 9s(1).$
Using Lemma \ref{lemma:5.8} gives  
\begin{equation}\label{eqn:5.11}
b(1) = \frac{(k-4)(29k-9)}{8(4k-1)(4k-9)}a(1).
\end{equation}
It is immediate that $b(1)\neq 0$ as $k$ is not 
a perfect square.

(2) Since
$\frac{m^2}{4k}=\frac{1}{16}$ and $\frac{9m^2}{4k}=\frac{9}{16}$ we have
$a(\frac{1}{16})=a(\frac{9}{16})=0.$ 
By Lemma \ref{lemma:5.8}, 
\[
b(\frac{1}{16}) = b(\frac{9}{16}) = 0.
\]
Next we assert that $b(1)=0.$ If $k=4$ this is immediate from 
(\ref{eqn:5.11}). If $k>4$ then $m>1.$ So $4m\leq k-1.$ Since 
$\frac{16m^2}{4k}=1$ we again have $a(1)=0$ and $b(1)=0$ by
using (\ref{eqn:5.11}). 

The proof of (3) is similar to that of (2).
\end{proof}

\begin{proposition}\label{proposition:5.11}
If $k$ is not a perfect square then
\[
\left\{
1,\omega,\omega^2,\dots,\omega^{k+3}, E,
\omega E,\omega^2 E
\right\}
\]
is a basis of $A(V_L^+)$. In particular, $\dim_\C A(V_L^+) = k+7$.
\end{proposition}
\begin{proof}
Note that ($B_2$) can be written as 
\begin{gather*}
(\omega-1)(\omega-\frac{1}{16})(\omega-\frac{9}{16}) J = 
(\omega-1)(\omega-\frac{1}{16})(\omega-\frac{9}{16})(4\omega^2-\omega).
\end{gather*}
From (\ref{ex1})-(\ref{ex2}) we see that
\[
b(\omega) J=r(\omega)a(\omega)+
\left\{
r(\omega)(\omega-4\omega^2)-p(\omega)
\right\}
s(\omega).
\]
Since $k$ is not a perfect square
 it follows from Lemma \ref{lemma:5.10} (1) that
$(x-1)(x-1/16)(x-9/16)$ and $b(x)$ are coprime. Then there exist
polynomials $\alpha(x)$ and $\beta(x)$ such that
\[
\alpha(x)(x-1)(x-\frac{1}{16})(x-\frac{9}{16})
+\beta(x)b(x) = 1.
\]
Thus
\begin{align*}
J = 
\alpha(\omega)&(\omega-1)(\omega-\frac{1}{16})
(\omega-\frac{9}{16})(4\omega^2-\omega)\\
&+\beta(\omega)\left\{
r(\omega)a(\omega)
+(r(\omega)(\omega-4\omega^2)-p(\omega))s(\omega)
\right\}.
\end{align*}
This shows that $A(V_L^+)$ is spanned by
$\left\{\omega^i,i\in\Z_{\geq 0},E,\omega E,\omega^2 E
\right\}$. Lemma \ref{lemma:5.6} then implies that
 $A(V_L^+)$ is spanned by
$\left\{1,\omega,\omega^2,\dots,\omega^{k+3},
E,\omega E,\omega^2 E
\right\}$ since $\deg\,\varphi = k+4$. Finally linear independence of
these vectors is clear because $A(V_L^+)$
has $k+7$ inequivalent irreducible modules which are the top levels
of the known irreducible modules for $V_L^+.$
\end{proof}

\begin{proposition}\label{proposition:5.12}
If $k=4m^2$ for positive integer $m$,
then the following set is a basis of $A(V_L^+)$;
\[
\left\{
1,\omega,\cdots,\omega^k,J,\omega J,\omega^2 J,
E,\omega E,\omega^2 E
\right\}.
\]
In particular, $\dim_\C A(V_L^+) =k+7$.
\end{proposition}
\begin{proof}

By Lemma \ref{lemma:5.10} (2), $b(x)$ has a factor $(x-1)(x-1/16)(x-9/16)$ and we can write
\[
b(x) = (x-1)(x-\frac{1}{16})(x-\frac{9}{16})c(x)
\]
where $c(x)$ is some polynomial. Even in this case, we still have
two relations;
\begin{gather*}
(\omega-1)(\omega-\frac{1}{16})(\omega-\frac{9}{16}) J = 
(\omega-1)(\omega-\frac{1}{16})(\omega-\frac{9}{16})(4\omega^2-\omega),\\
(\omega-1)(\omega-\frac{1}{16})(\omega-\frac{9}{16})c(\omega)J
=r(\omega)a(\omega)+
\left\{
r(\omega)(\omega-4\omega^2)-p(\omega)
\right\}
s(\omega)
\end{gather*}
(see ($B_2$) and (\ref{ex1})-(\ref{ex2})).
By eliminating $J$ we obtain
\[
b(\omega)(4\omega^2 -\omega)
= r(\omega)a(\omega)+
\left\{
r(\omega)(\omega-4\omega^2)-p(\omega)
\right\}
s(\omega).
\]
Using the definition of $b(x)$ we see that
\begin{align*}
\{r(\omega)&-4\omega^4 +\omega\}a(\omega)\\
&+\left\{
r(\omega)(\omega-4\omega^2)-p(\omega)
-(q(\omega)-r(\omega)+\omega-4\omega^2)(4\omega^2-\omega)\right\}
s(\omega)= 0
\end{align*}
or that 
\[
(r(\omega)-4\omega^2+\omega)a(\omega)
+
\left\{
-p(\omega)-q(\omega)(4\omega^2-\omega)+(4\omega^2-\omega)^2
\right\}
s(\omega) = 0.
\]

A direct calculation gives 
\[
-p(x)-q(x)(4x^2-x)+(4x^2-x)^2 = 0
\]
and 
\begin{equation}\label{eqn:5.12}
r(x) -4x^2+x = \frac{9(4x-k)((32k-12)x-3k+3)}{8(4k-1)(4k-9)}.
\end{equation}
Thus
\begin{equation}\label{eqn:5.12'}
(\omega-\frac{k}{4})(\omega-\frac{3k-3}{32k-12})a(\omega) = 0.
\end{equation}

Recall the definition of $\varphi(x)$ from (\ref{ex3}). 
Suppose that $\frac{3k-3}{32k-12}$ is not a root of $\varphi(x).$
Then $x-\frac{3k-3}{32k-12}$ and $\varphi(x)$ are relatively
prime and 
there exist polynomials $f(x)$ and $g(x)$ such that  
$(x-\frac{3k-3}{32k-12})f(x)+\varphi(x)g(x)=1.$ So 
$$(x-\frac{k}{4})a(x)=(x-\frac{k}{4})(x-\frac{3k-3}{32k-12})a(x)f(x)+(x-\frac{k}{4})a(x)\varphi(x)g(x).$$
By Lemma \ref{lemma:5.6}
\[
(\omega-\frac{k}{4})a(\omega) = 0.
\]
So $\o^{k+1}$ is a linear combinations of $\o^i$ for $i\leq k.$ 
We then use the relations ($B_1$), ($B_2$),($L_1$) and ($L_2$) to conclude
that 
\[
\left\{
1,\omega,\cdots,\omega^k,J,\omega J,\omega^2 J,
E, \omega E, \omega^2 E
\right\}
\]
is a spanning set of $A(V_L^+)$. Since $A(V_L^+)$ has $k+7$ known inequivalent
irreducible modules already we see immediately that this
spanning set is a basis.

It remains to prove that $\frac{3k-3}{32k-12}$ is not a root of
$\varphi(x).$ Note that the roots of $\varphi(x)$ are
$$1, \frac{1}{16},\frac{9}{16},\frac{k}{4}, \frac{i^2}{4k}, i=0,...,k-1.$$
It is easy to see that $\frac{3k-3}{32k-12}\ne 0,1,\frac{1}{16},
\frac{9}{16}, \frac{k}{4}.$ Suppose that 
\[
\frac{3k-3}{32k-12}=\frac{i^2}{4k}
\]
or
\[
\frac{12m^2-3}{128m^2-12}=\frac{i^2}{16m^2}
\]
for some $ \geq 1\leq k-1.$ 
Then
\[
i^2 = \frac{12m^2(4m^2-1)}{32m^2-3}.
\]
Let $d$ be the greatest common divisor of $12m^2(4m^2-1)$ and $32m^2-3$,
then $d$ divides $-24m^2(4m^2-1)+3m^2(32m^2-3)$ or $15m^2.$ 
Since $\frac{12m^2(4m^2-1)}{32m^2-3}$ is an integer we see that 
$d=32m^2-3$ and $32m^2-3$ divides $15m^2.$ This is impossible for any positive integer $m$. 
So we have a contradiction.
\end{proof}

\begin{proposition}\label{proposition:5.13}
If $k=(2m+1)^2$ for positive integer $m$,
then the following set is a basis of $A(V_L^+)$;
\[
\left\{
1,\omega,\cdots,\omega^{k+2},J, E,\omega E,\omega^2 E
\right\}.
\]
In particular, $\dim_\C A(V_L^+) = k+7$.
\end{proposition}
\begin{proof}

By Lemma \ref{lemma:5.8} (3) we can write $b(x)=(x-1)c(x)$
such that the polynomials $c(x)$ and $(x-1/16)(x-9/16)$ are coprime. 
Use the following relations
\begin{gather}
(\omega-1)(\omega-\frac{1}{16})(\omega-\frac{9}{16}) J = 
(\omega-1)(\omega-\frac{1}{16})(\omega-\frac{9}{16})(4\omega^2-\omega),\notag\\
(\omega-1)c(\omega)J=r(\omega)a(\omega)+
\left\{
r(\omega)(\omega-4\omega^2)-p(\omega)
\right\}
s(\omega)\label{ex4}
\end{gather}
(see ($B_2$) and (\ref{ex1})-(\ref{ex2})) to eliminate $J$ and
to obtain 
$$(\omega-\frac{1}{16})(\omega-\frac{9}{16})(\omega-\frac{k}{4})(\omega-\frac{3k-3}{32k-12})a(\omega) = 0.$$
Again one can show that $\frac{3k-3}{32k-12}$ is not a root of $\varphi(x)$
as in the proof of Proposition \ref{proposition:5.12}. Thus we have
$$(\omega-\frac{1}{16})(\omega-\frac{9}{16})(\omega-\frac{k}{4})a(\omega) = 0.$$
and $\o^{k+3}$ is a linear combination of $\o^i$'s for $i=0,...,k+2.$

Let the polynomials $\a(x)$ and $\b(x)$ satisfy 
$$\a(x)c(x)+\b(x)(x-1/16)(x-9/16)=1.$$
Combining this with relations (\ref{ex4}) gives
\begin{align*}
(\o-1)J=\b(\o)&(\o-1) (\o-\frac{1}{16})(\o-\frac{9}{16})(4\o^2-\o)\\ 
&+\a(\omega)\left\{
r(\omega)a(\omega)
+(r(\omega)(\omega-4\omega^2)-p(\omega))s(\omega)
\right\}.
\end{align*}
Thus $A(V_L^+)$ is spanned by
$\left\{\omega^i,i=0,...,k+2,J, E,\omega E,\o^2 E\right\}$. Again
the known simple modules for $A(V_L^+)$ implies that
this in fact is a basis of $A(V_L^+).$
\end{proof}

\begin{remark} We can determine all relations in $A(V_L^+).$ From the proofs
of Propositions \ref{proposition:5.11}-\ref{proposition:5.13} it is enough
to give the exact expression of $s(x)$ in (\ref{eqn:5.10}). Since
$s(x)$ is a polynomial of degree less than or equal to $2$ there
are at most $3$ coefficients to be determined. This can be done by
evaluating (\ref{eqn:5.10}) on the top levels of $V_L^-,$ $V_L^{T_{1},+}$
and $V_L^{T_{1},-}.$ We leave the detail to the reader.
\end{remark}

\subsection{Classification of irreducible modules}

Recall the known irreducible modules for $V_L^+$ from Section 2. 
We finally have the following classification result:
\begin{theorem} Let $L=\Z\a$ be a even positive definite lattice of rank
1 such that $\a$ has square length $2k.$ Then 
$$\{V_{L}^{\pm}, V_{L+\frac{\a}{2}}^{\pm}, V_L^{T_i,\pm},V_{L+\frac{r}{2k}\alpha}|i=1,2, r=1,...,k-1\}$$
gives a complete irreducible modules for $V_L^+.$ Moreover,
any admissible irreducible $V_L^+$-module is an ordinary module.
\end{theorem}

\begin{proof} If $k\ne 1$ then $A(V_L^+)$ is a commutative algebra
of dimension $k+7$ (see Section 5). So $A(V_L^+)$ has at most 
$k+7$ simple modules. Since $A(V_L^+)$ has $k+7$ known
inequivalent simple modules already we conclude that $A(V_L^+)$
is a semisimple algebra of dimension $k+7.$ Using the one to one correspondence
result (see Theorem \ref{theorem:2.3}) we see that
$V_L^+$ has exactly $k+7$ inequivalent irreducible admissible 
modules which are ordinary modules.

If $k=1,$ $V_L^+$ is isomorphic to the the lattice vertex operator
algebra $V_{L'}$ where $L'$ is a rank one positive definite lattice
spanned by $\beta$ whose square length is 8 (see [DG]). Then it
follows from a result in [D1] that $V_L^+$ has exactly 8 irreducible
modules. Since $V_{L'}$ is rational [DLM1] every irreducible admissible
module is an ordinary module. 
\end{proof}

\begin{remark} If $k=2$ the vertex operator algebra $V_L^+$ is isomorphic
to $L(1/2,0)\otimes L(1/2,0)$ (see Lemma 3.1 of [DGH]) where
$L(1/2,h)$ is the irreducible highest weight module for the
Virasoro algebra with central charge $1/2$ and highest weight $h.$
So the classification of irreducible modules in this case
also follows from the classification of irreducible modules
for the vertex operator algebra $L(1/2,0).$ One can easily see that 
all the irreducible modules for $L(1/2,0)\otimes L(1/2,0)$
are $L(1/2,h_1)\otimes L(1/2,h_2)$ where $h_i=0,1/16, 1/2.$ One can find 
in [DGH] the identification of these modules with the  modules for
$V_L^+$ listed in the theorem.
\end{remark}

\end{document}